\def\of#1{\left(#1\right)}
\def\ofc#1{\left\{#1\right\}}
\def\ofb#1{\left[#1\right]}

\def\Devaney{1}
\def\Debski{2}
\def\Fokkink{3}
\def\FortMcCord{4}
\def\KeeslingOne{5}
\def\KeeslingTwo{6}
\def\Kennedy{7}
\def\treelike{8}
\def\solenoids{9}
\def\RogersTollefson{10}
\def\RogersOne{11}
\def\RogersTwo{12}
\def\Scheffer{13}
\def\Watkins{14}

\input epsf
\input amstex
\documentstyle{amsppt}

\topmatter

\title Homotopy classes of maps between Knaster continua \\
\endtitle

\rightheadtext {Homotopy classes of maps between Knaster continua }

\author Piotr Minc \endauthor

\address Department of Mathematics, Auburn University,
Auburn, Alabama
36849\endaddress
\email mincpio\@mail.auburn.edu \endemail
\thanks This research was supported in part by NSF grant
\# DMS-9505108.
\endthanks

\abstract
By a Knaster continuum we understand the inverse limit of
copies of $\ofb{0,1}$ with open bonding maps.
We prove that
for any two Knaster continua $K_1$ and $K_2$,
there are $2^{\aleph_0}$ distinct homotopy types of maps
of $K_1$ onto $K_2$
that map the endpoint of $K_1$ to the endpoint of $K_2$.
\endabstract
\endtopmatter

\document


\head 1. Introduction \endhead

Let $\Bbb R$ denote the set of real numbers
and let $I$ denote the interval $\ofb{0,1}$. For any real number $t$,
let $\ofb{t}$ denote the greatest integer less than or equal to $t$.
Let
$v:\Bbb R\to I$ be defined by the formula
$$v\of{t}=\cases
t-\ofb{t},&\text{if $\ofb{t}$ is even}\\
\ofb{t}+1-t,&\text{if $\ofb{t}$ is odd}.
\endcases
$$
For each positive integer $n$, let $g_n:I\to I$ be defined by the
formula $g_n\of{t}=v\of{nt}$. Observe that $g_n$ stretches $n$ times
and then folds the resulting interval $\ofb{0,n}$ onto $\ofb{0,1}$.
The map
$g_2$ is the very well
known \lq\lq roof-top\rq\rq{} map. For any two positive integers
$m$ and $n$, $g_m\circ g_n=g_{mn}$. Consequently, $g_n$ and
$g_m$ commute
(see, for example, \cite{\treelike{\rm , Proposition 2.2}}).

Let $N=\ofc{n_1,n_2,\dots}$ be a sequence of integers $>1$.
Consider the inverse sequence
$$
I@<g_{n_1}<<I@<g_{n_2}<<I@<g_{n_3}<<I@<g_{n_4}<<\dots\tag{*}
$$
By the Knaster continuum associated with the sequence $N$
we will understand the inverse limit of \thetag{*}.
Observe that the same Knaster continuum can be associated with
two different sequences. For example the inverse limit
does not change if we replace $N$ by
$\ofc{n_1 n_2\dots n_{j_1}, n_{j_1+1}\dots n_{j_2},\dots}$.
However, it should be noted here that there are
$2^{\aleph_0}$ topologically distinct
Knaster continua \cite{\Debski}(see  also \cite{\Watkins}).

For a
Knaster continuum $K$, let $e$ denote the endpoint
$\of{0,0,\dots}$. By $\pi_i$ we will understand the projection
of $K$ onto the $i$-th component in the inverse system
($i=0,1,\dots$).

Let $S^1$ denote the unit circle in the complex plane and
let $\tilde g_n:S^1\to S^1$ be defined by
$\tilde g_n\of{z}=z^n$. By the solenoid
associated with the sequence $N$
we will understand the inverse limit of
$
S^1@<\tilde g_{n_1}<<
S^1@<\tilde g_{n_2}<<
S^1@<\tilde g_{n_3}<<
S^1@<\tilde g_{n_4}<<\dots
$.
We will use $e$ to denote $\of{1,1,\dots}$.

Note that the same letter $e$ is used here to denote different
objects:
$\of{0,0,\dots}$ in different Knaster continua and
$\of{1,1,\dots}$ in different solenoids.
For instance,
if $f$ is a map between two Knaster continua $K_1$ and $K_2$,
then the equality $f\of{e}=e$ means that $f$ maps the
endpoint $\of{0,0,\dots}$ of $K_1$ to the
endpoint $\of{0,0,\dots}$ of $K_2$.

Knaster continua and solenoids are the simplest examples
of indecomposable continua. They appear naturally as
attractors of dynamical systems (see for example \cite{\Devaney}
and \cite{\Kennedy}).
The goal of this paper is to study homotopy classes of maps between
Knaster continua.
It follows from \cite{\RogersTwo} (see also \cite{\FortMcCord}
and \cite{\RogersTollefson}) that
each map between two Knaster continua is homotopic to a map induced from
a commutative diagram. Hopotopy classes of maps between
Knaster continua coincide, therefore, with those of the induced
map. The easy examples of commutative diagrams could be
constructed by using the commutativity of $g_n$ and $g_m$.
We will first consider such examples.

Let $M=\ofc{m_1,m_2,\dots}$ be another sequence of positive integers.
Suppose $\ofc{k_0,k_1,k_2,\dots}$ is the third
sequence of positive integers such that
$k_{i-1}n_i=m_ik_i$ for each positive integer $i$.
Then the following
diagram commutes.

$$
\CD
I@<g_{n_1}<<I@<g_{n_2}<<I@<g_{n_3}<<I@<g_{n_4}<<\dots \\
@Vg_{k_0}VV @Vg_{k_1}VV @Vg_{k_2}VV @Vg_{k_3}VV\\
I@<<g_{m_1}<I@<<g_{m_2}<I@<<g_{m_3}<I@<<g_{m_4}<\dots
\endCD
\tag{**}
$$
Let $K_1$ and $K_2$ be Knaster continua associated with
the sequences $N$ and $M$, respectively. Observe that
the diagram \thetag{**} induces a continuous map
$f:K_1\to K_2$. Note that if $x=\of{x_0,x_1,x_2,\dots}$
then
$f\of{x}=\of{g_{k_0}\of{x_0}, g_{k_1}\of{x_1},g_{k_2}\of{x_2},
\dots}$.

We will say that a map $g:K_1\to K_2$ is
{\it naturally induced\/} if there are sequences
$\ofc{j_0,j_1,j_2,\dots}$ and $\ofc{i_0,i_1,i_2,\dots}$
such that $0\le j_0<j_1<j_2<\dots$ and
$$g\of{x}=\of{g_{i_0}\of{x_{j_0}}, g_{i_1}\of{x_{j_1}},
g_{i_2}\of{x_{j_2}},
\dots}.$$
The map induced by the diagram \thetag{**} is an example
of a naturally induced map. In the general case, the
vertical arrows in the diagram \thetag{**} may be replaced by
diagonal ones.
The following proposition is a simple consequence of the definition.
\proclaim{Proposition 1.1} Suppose $K_1$ and $K_2$ are
Knaster continua associated with the sequences
$\ofc{n_1,n_2,\dots}$ and
$\ofc{m_1,m_2,\dots}$, respectively.
Let $\ofc{i_0,i_1,\dots}$
and $\ofc{j_0,j_1,\dots}$ be two sequences
of non negative
integers
with $j_0<j_1<j_2<\dots$.
Let
$g:K_1\to K_2$ be the naturally induced map defined by
$$g\of{x}=\of{g_{i_0}\of{x_{j_0}}, g_{i_1}\of{x_{j_1}},
g_{i_2}\of{x_{j_2}},
\dots}$$
for each $x=\of{x_0,x_1,x_2,\dots}\in K_1$.
Then
$${i_k}=
\frac{{i_0}n_{j_0+1}n_{j_0+2}\dots n_{j_k}}{m_1 m_2 \dots m_k}
$$
for each positive integer $k$.
\endproclaim

\proclaim{Corollary 1.2}
Suppose $K_1$ and $K_2$ are two Knaster continua and
$\pi^{\prime\prime}_0$ denote the projection
of $K_2$ onto the $0-th$ factor in the
inverse sequence defining $K_2$.
Let $g$ and $f$ be two naturally
induced maps between of $K_1$ into $K_2$
such that $\pi^{\prime\prime}_0\circ g=\pi^{\prime\prime}_0\circ f$.
Then $g=f$.
\endproclaim

\proclaim{Corollary 1.3} The set of naturally
induced maps between two Knaster continua is countable.
\endproclaim

The notion of naturally induced maps may be introduced for solenoids.
Suppose $\Sigma_1$ and $\Sigma_2$ are solenoids  associated with
the sequences $N$ and $M$, respectively.
We will say that a map $\tilde g:\Sigma_1\to \Sigma_2$ is
{\it naturally induced\/} if there are sequences
$\ofc{j_0,j_1,j_2,\dots}$ and $\ofc{i_0,i_1,i_2,\dots}$
such that $0\le j_0<j_1<j_2<\dots$ and
$$\tilde g\of{x}=\of{\tilde g_{i_0}\of{x_{j_0}},
\tilde g_{i_1}\of{x_{j_1}},
\tilde g_{i_2}\of{x_{j_2}},
\dots}.$$
The statements corresponding to 1.1, 1.2 and 1.3 are
true for solenoids.
Observe that any naturally induced map, either between
two Knaster continua or between two solenoids, maps
$e$ to $e$. Suppose $f$ is a map with the property
that $f\of{e}=e$. Must $f$ be homotopic to
a naturally induced map? In case of solenoids, the
answer is positive. The following proposition
follows from \cite{\Scheffer}, \cite{\KeeslingOne} and
\cite{\KeeslingTwo}
(see also \cite{\solenoids {\rm, Proposition 3}}).

\proclaim{Proposition 1.4}
Suppose $f$ is a map of a solenoid $\Sigma_1$ into a solenoid
$\Sigma_2$ such that
that $f\of{e}=e$. Then $f$ is homotopic to
a naturally induced map.
\endproclaim

A similar statement for Knaster continua
would be false.
For instance, any Knaster
continuum $K_1$ can be always mapped onto any Knaster
continuum $K_2$
(see \cite{\RogersOne}).
On the other hand, there is no naturally
induced maps of $K_1$ of $K_2$
if the sequences defining $K_1$ and
$K_2$ have different prime factors (see Proposition 1.1).

It follows from Proposition 1.4 that there is only
countably many homotopy classes of maps between solenoids
mapping $e$ into $e$. One could expect the
corresponding theorem for Knaster continua since the
mapping structure of solenoids is usually richer.
After all, any solenoid can be mapped onto any Knaster
continuum, but there is no continuous map of a
chainable continuum onto a solenoid.
We will show, however, that in case of Knaster continua
there are uncountably many homotopy types of maps
mapping the endpoint of the domain onto the endpoint
of the range. More precisely,
we will prove the following theorem.
\proclaim{Theorem 1.5}
For any Knaster continua $K_1$ and $K_2$,
there are $2^{\aleph_0}$ distinct homotopy types of maps
$f:K_1\to K_2$
such that $f\of{e}=e$.
\endproclaim
It should be noted here that the maps in the above theorem cannot
be replaced by homeomorphisms even in the case where $K_1=K_2$.
It follows from \cite{\Fokkink{\rm, Theorem 3.3, p. 53}} that there is only
countably many homotopy types of homeomorphisms of
any Knaster continuum onto itself.

\head 2. Constructing Maps between Knaster continua.\endhead
\midinsert
$$\vbox{\epsffile{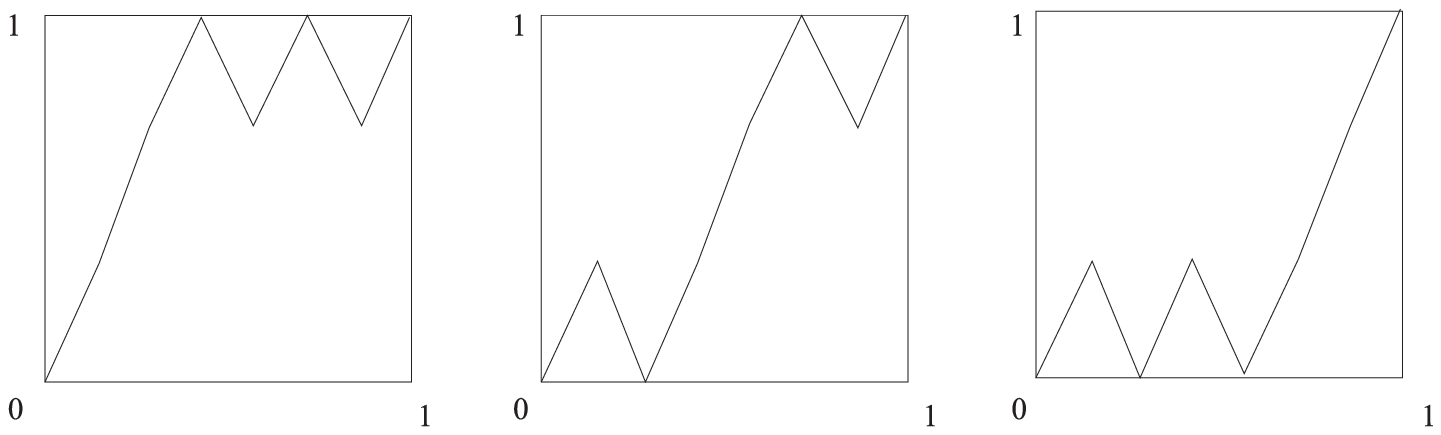}}$$
\centerline{Figure}
\endinsert
We will prove Theorem 1.5 by defining $2^{\aleph_0}$
maps of $K_1$ onto $K_2$ induced by commutative diagrams.
We will construct the vertical arrows $f_0, f_1, f_2,\dots$
one by one. We will start
from $f_0$ equal to the identity on $I$. We
will then use Lemma 2.1 repeatedly to define $f_1, f_2,\dots$
in many different ways. The figure depicts the graphs
of three maps $f_1$ such that
$f_0\circ g_7=g_3\circ f_1$. Similar maps will be used in
Lemma 2.1.

\proclaim{Lemma 2.1} Let $m$, $n$ and $q$ be positive
integers such that $\of{m+2}q\le n$. Suppose $f_0$ is a map of $I$
onto itself.
Then for each
$i=0,1,\dots,q-1$ there is a map $f_1:I\to I$ such that
\roster
\item{} $f_1\of{0}=0$ if $f_0\of{0}=0$,
\item{} $f_0\circ g_n=g_m\circ f_1$,
\item{} $f_1\of{\ofb{0,\frac{i}{q}}}\subset\ofb{0,\frac{1}{m}}$,
\item{} $f_1\of{\ofb{\frac{i}{q},\frac{i+1}{q}}}=I$ and
\item{} $f_1\of{\ofb{\frac{i+1}{q},1}}\subset\ofb{\frac{m-1}{m},1}$.
\endroster
\endproclaim
\demo{Proof of {\rm 2.1}}
Let $k$ be the least non negative integer such that
$\frac{k}{n}\ge\frac{i}{q}$. Since $\of{m+2}q\le n$, it follows
that $\frac{k+m+1}{n}<\frac{i+1}{q}$ and consequently
$$\ofb{\frac{k}{n},\frac{k+m+1}{n}}\subset
\ofb{\frac{i}{q},\frac{i+1}{q}}.$$

Let $a,b\in I$ be such that
$f_0\of{a}=0$ and $f_0\of{b}=1$.
Since $g_n$ maps each of the intervals
$\ofb{\frac{k+j}{n},\frac{k+j+1}{n}}$ onto $I$,
for each $j=0,\dots,m$,  there is a point
$t_j\in\ofb{\frac{k+j}{n},\frac{k+j+1}{n}}$ such that
$g_n\of{t_j}=a$ if $j$ is even and
$g_n\of{t_j}=b$ if $j$ is odd. Define $f_1$ by
$$
f_1\of{t}=\frac{1}{m}
\cases
f_0\circ g_n\of{t},&\text{for $0\le t\le t_1$}\\
j+1-f_0\circ g_n\of{t},
&\text{for an odd $j=1,\dots,m-1$ and $t_j\le t\le t_{j+1}$ }\\
j+f_0\circ g_n\of{t},
&\text{for an even $j=2,\dots,m-1$ and $t_j\le t\le t_{j+1}$}\\
m-1 +f_0\circ g_n\of{t},&\text{for $t_m\le t\le1$ if $m$ is odd}\\
m-f_0\circ g_n\of{t},&\text{for $t_m\le t\le1$ if $m$ is even}.
\endcases
$$
One can verify that so defined $f_1$ has the required properties.
(Use the equality $g_m\of{x}=v\of{mx}$ to verify 2.)
\enddemo

\definition{Construction of $f_t$}
Let $K_1$ and $K_2$ be the Knaster continua associated with sequences
$\ofc{n_1,n_2,\dots}$ and $\ofc{m_1,m_2,\dots}$, respectively.
Observe that $K_1$ does not change if the sequence
$\ofc{n_1,n_2,\dots}$ is replaced by
a sequence of products of finite blocks of consecutive elements of
its elements.
So we may assume that
$$\of{m_j+2}j<n_j\quad \text{ for each positive integer $j$.}$$
Let $t\in I$.
For each positive integer $j$, let $i\ofb{t,j}$ be an integer such
that $0\le i\ofb{t,j}<j$ and
$t\in\ofb{\frac{i\ofb{t,j}}{j},\frac{i\ofb{t,j}+1}{j}}$.
We will define a sequence $f^t_0, f^t_1, f^t_2,\dots$ of maps
of $I$ onto itself such that
\roster
\item{}$f^t_{j}\of{0}=0$ ,
\item{}$f^t_{j-1}\circ g_{n_j}=g_{m_j}\circ f^t_{j}$,
\item{} $f^t_{j}\of{\ofb{0,\frac{i\ofb{t,j}}{j}}}\subset
\ofb{0,\frac{1}{m_j}}$,
\item{}$f^t_{j}\of{\ofb{\frac{i\ofb{t,j}}{j},\frac{i\ofb{t,j}+1}{j}}}=I$ and
\item{} $f^t_{j}\of{\ofb{\frac{i\ofb{t,j}+1}{j},1}}\subset
\ofb{\frac{m_j-1}{m_j},1}$.

\endroster
for each positive integer $j$.

Let $f^t_0$ be the identity on $I$. Suppose
$f^t_0,\dots,f^t_{j-1}$ have been defined. Use Lemma 2.1 with
$m=m_j$, $n=n_j$, $q=j$,
$f_0=f^t_{j-1}$ and $i=i\ofb{t,j}$. Set
$f^t_{j}$ to be $f_1$ obtained from the lemma.
Observe that conditions
1-5 follow from the corresponding conditions in the lemma.

Let $f^t:K_1\to K_2$ be the function induced by the sequence
$f^t_0, f^t_1, f^t_2,\dots$.
\enddefinition
The following proposition is a simple consequence of the
construction.
\proclaim{Proposition 2.2} For each $t\in I$, $f^t$ is a
continuous map of $K_1$ onto $K_2$ such that $f^t\of{e}=e$.
\endproclaim
\proclaim{Proposition 2.3} Suppose $t,s\in I$
and $t\ne s$. Then $f^t$ is not homotopic to $f^s$.
\endproclaim
\demo{Proof of {\rm 2.3}}
Suppose $f^t$ is homotopic to $f^s$ for some $t,s\in I$
such that $t<s$. Let
$H:K_1\times I\to K_2$ be the homotopy between
$f^t$ and $f^s$. Consider the homotopy
$h=\pi^{\prime\prime}_0\circ H:K_1\times I\to I$,
where $\pi^{\prime\prime}_0$ denote
the projection of $K_2$ onto the $0-th$ factor in the
inverse sequence defining $K_2$. By compactness, there
is a sequence of numbers $z_0=0<z_1<\dots<z_{\ell-1}<z_\ell=1$
such that the set $h\of{\ofc{x}\times\ofb{z_{k-1},z_k}}$
does not contain $I$ for each $x\in K_1$ and $k=1,\dots,\ell$.

Let $j$ be an integer such that
$$m_1 m_2\dots m_{j-1}>\ell \quad\text{and}\quad
\frac{3}{j}<s-t.$$
Since $2j<n_j$, there is an integer $q$ such that
$$\frac{i\ofb{t,j}+1}{j}\le\frac{2q}{n_j}\le\frac{i\ofb{t,j}+2}{j}.$$
By condition 5 of the construction of $f^t$, we have that
$$f^t_{j}\of{\frac{2q}{n_j}}\in
\ofb{\frac{m_j-1}{m_j},1}.\tag{i}$$
Since
$\frac{3}{j}<s-t$ and
$t\in\ofb{\frac{i\ofb{t,j}}{j},\frac{i\ofb{t,j}+1}{j}}$,
we have that
$\frac{i\ofb{t,j}+3}{j}<s$.
Since $s\in\ofb{\frac{i\ofb{s,j}}{j},\frac{i\ofb{s,j}+1}{j}}$,
it follows that
$\frac{i\ofb{t,j}+2}{j}<\frac{i\ofb{s,j}}{j}$ and
consequently
$\frac{2q}{n_j}\in\ofb{0,\frac{i\ofb{s,j}}{j}}$.
By condition 3 of the construction of $f^s$, we have that
$$f^s_{j}\of{\frac{2q}{n_j}}\in
\ofb{0,\frac{1}{m_j}}.\tag{ii}$$
Since $g_{n_j}\of{\frac{2q}{n_j}}=0$ and $f^s_{j-1}\of{0}=0$,
condition 2 of the construction of $f^s$ implies that
$g_{m_j}\circ f^s_{j}\of{\frac{2q}{n_j}}=0$. It follows from
\thetag{ii} that
$$f^s_{j}\of{\frac{2q}{n_j}}=0
.\tag{iii}
$$

Let $\pi^\prime_j$
denote
the projection of $K_1$ onto the $j-th$ factor in the
inverse sequence defining $K_1$. Similarly, let
$\pi^{\prime\prime}_j$
denote
the projection of $K_2$ onto the $j-th$ factor in the
inverse sequence defining $K_2$.

Take a point $y\in K_1$ such that
$\pi^\prime_j\of{y}=\frac{2q}{n_j}$.
It follows from \thetag{i} and \thetag{iii} that
$$\pi^{\prime\prime}_j\circ f^s\of{y}=0
\quad\text{ and }\quad
\pi^{\prime\prime}_j\circ f^t\of{y}\in
\ofb{\frac{m_j-1}{m_j},1}.
\tag{iv}$$

Let $p$ denote the product $m_1 m_2\dots m_{j-1}$ and
let $r=pm_j$. By \thetag{iv}, there are numbers
$w_0, w_1,\dots w_p\in I$ such that
$0=w_0<w_1<\dots<w_p$ and
$$\pi^{\prime\prime}_j\circ H\of{y,w_\lambda}=\frac{\lambda}{r}
\quad\text{for
each $\lambda=0,1,\dots,p$.}$$
Observe that
$h\of{y,w_\lambda}=\pi^{\prime\prime}_0\circ H\of{y,w_\lambda}=
g_r\circ \pi^{\prime\prime}_0\circ H\of{y,w_\lambda}=
g_r\of{\frac{\lambda}{r}}$. Thus,
$$\ofc{h\of{y,w_{\lambda-1}},h\of{y,w_\lambda}}=\ofc{0,1}
\quad\text{for
each $\lambda=1,\dots,p$}.\tag{v}$$
Recall that $p>\ell$ by the choice of
$j$. Hence, there are integers $k$ and $\lambda$ such that
$1\le k\le\ell$, $1\le\lambda\le p$ and
$$w_{\lambda-1},w_{\lambda}\in\ofb{z_{k-1},z_k}.$$
Now, \thetag{v}
contradicts the choice of $z_0,z_1,\dots,z_{\ell}$.

\enddemo

\demo{Proof of Theorem {\rm 1.5}} The theorem follows from
Propositions 2.2 and 2.3.
\enddemo



\enddocument

\bye